\documentclass[11pt,a4paper]{amsart}
\usepackage[T1]{fontenc}
\usepackage{lmodern}
\usepackage{amsmath,amssymb,mathtools}
\usepackage{microtype,enumitem}
\usepackage[a4paper,textwidth=148mm,textheight=225mm,centering,headheight=14pt]{geometry}
\usepackage[hidelinks]{hyperref}
\hypersetup{pdftitle={A numerical fixed-point criterion for automorphism groups of curves, with an application to the generalized Suzuki curve},
 pdfauthor={Saeed Tafazolian},
 pdfkeywords={Algebraic curves; automorphism groups; p-rank; ramification groups; generalized Suzuki curve}}
\numberwithin{equation}{section}
\newtheorem{theorem}{Theorem}[section]
\newtheorem{proposition}[theorem]{Proposition}
\newtheorem{lemma}[theorem]{Lemma}

\newtheorem*{theoremA}{Theorem A}
\newtheorem*{theoremB}{Theorem B}
\newtheorem*{theoremC}{Theorem C}
\newtheorem*{corollaryD}{Corollary D}
\theoremstyle{remark}
\newtheorem{remark}[theorem]{Remark}
\newtheorem{example}[theorem]{Example}
\newcommand{\Aut}{\operatorname{Aut}}
\newcommand{\Syl}{\operatorname{Syl}}
\newcommand{\PGU}{\operatorname{PGU}}
\newcommand{\PGL}{\operatorname{PGL}}
\newcommand{\Tr}{\operatorname{Tr}}
\newcommand{\tr}{\operatorname{tr}}
\newcommand{\F}{\mathbb F}
\newcommand{\PP}{\mathbb P}
\newcommand{\cX}{\mathcal X}
\newcommand{\GS}{\mathrm{GS}}

\title[A numerical fixed-point criterion]
{A numerical fixed-point criterion and a characterization of the generalized Suzuki curve}
\author{Saeed Tafazolian and Yuri da Silva}
\address{IMECC, Universidade Estadual de Campinas (UNICAMP),
Campinas, SP, Brazil}
\email{saeed@unicamp.br}
\email{y225979@dac.unicamp.br}
\date{}

\begin{document}
\begin{abstract}
Let $G$ be a finite automorphism group of a curve of genus $g\ge2$
over an algebraically closed field of odd characteristic $p$, and let
$S$ be a Sylow $p$-subgroup of $G$, of order $s>1$. We show that if
$s>p(g-1)/(p-2)$ and no divisor $d>1$ of the index $[G:S]$ satisfies
$d\equiv1\pmod s$ (for instance, if $|G|<s(s+1)$), then the curve has
$p$-rank zero, $S$ is normal in $G$, $G=S\rtimes C$ with $C$ cyclic,
and $G$ fixes exactly one point. The proof combines Nakajima's bound
with a counting argument and uses neither solvability nor
classification results. As an
application, let $q_0=p^\ell$ and $q=p^{2\ell-1}$ with $\ell\ge2$. We
prove that a curve of genus $q_0(q-1)/2$ admits an automorphism group
of order $q^2(q-1)$ if and only if it is isomorphic to the generalized
Suzuki curve $y^q-y=x^{q_0}(x^q-x)$. Over $\F_q$ the same holds up to
$\F_q$-isomorphism when the automorphisms are defined over $\F_q$, and
twists show that this condition cannot be dropped. This sharpens an
earlier characterization, which assumed in addition that the
automorphism group fixes an $\F_q$-rational point.
\end{abstract}
\maketitle

\section{Introduction}\label{sec:intro}

Let $K$ be an algebraically closed field of characteristic $p>2$,
and let $\cX/K$ be a smooth projective curve of genus $g\ge2$, with
$p$-rank $\gamma(\cX)$. A finite group $G\le\Aut(\cX)$ whose order
is divisible by $p$ need not fix a point of $\cX$. On the other hand,
the stabilizer $G_P$ of a point $P$ has a unique Sylow $p$-subgroup,
and large $p$-subgroups of point stabilizers are strongly restricted.
Stichtenoth \cite{StArch} bounded their order in terms of $g$, and
Giulietti and Korchm\'aros \cite{GK10} proved that if $\Aut(\cX)$
fixes no point of $\cX$ while some point stabilizer contains a
$p$-subgroup of order greater than $pg/(p-1)$, then $\cX$ is one of a
few classical curves. Their proof relies on the classification of
doubly transitive permutation groups with cyclic two-point
stabilizers \cite{KOS}. On the side of the $p$-rank, Nakajima
\cite{Nak} proved that every $p$-subgroup $S\le\Aut(\cX)$ satisfies
$|S|\le p(\gamma(\cX)-1)/(p-2)$ whenever $\gamma(\cX)\ge2$. For large
groups, Giulietti and Korchm\'aros \cite{GK19} showed that if $G$ is
solvable and $|G|>252g^2$, then $\gamma(\cX)=0$ and $G$ fixes a point
of $\cX$; see \cite{GKT} for recent bounds on automorphism groups of
curves of positive $p$-rank in terms of $\gamma(\cX)$.

In this note we combine Nakajima's bound with an elementary counting
argument to obtain a numerical criterion which forces a finite group
of automorphisms to fix a point. It involves only the order of a
Sylow $p$-subgroup and the order of the group, and it uses neither
solvability nor classification results. In contrast with
\cite{GK19}, the size hypothesis concerns the Sylow $p$-subgroup
rather than the whole group. This is the form needed in our
application, where nothing but the order of the group is known and
$|G|\le 2g^2$ (Remark~\ref{rem:GK19}).

\begin{theoremA}
Let $G\le\Aut(\cX)$, let $S\in\Syl_p(G)$, and put $s=|S|>1$. Assume
that
\begin{enumerate}[label=\textup{(\roman*)},leftmargin=*]
\item $s>\dfrac{p}{p-2}\,(g-1)$, and
\item no divisor $d>1$ of $[G:S]$ satisfies $d\equiv1\pmod s$.
\end{enumerate}
Then $\gamma(\cX)=0$, $S$ is normal in $G$, $G=S\rtimes C$ with $C$
cyclic of order $[G:S]$, and $G$ has exactly one fixed point on
$\cX$. If, moreover, $K=\overline{\F}$ for a finite field $\F$, the
curve $\cX$ is defined over $\F$, and $G$ is stable under conjugation
by the $\F$-Frobenius, then this point is $\F$-rational.
\end{theoremA}

Condition~(ii) depends only on $s$ and $|G|$. It holds in particular
whenever $|G|<s(s+1)$, since a divisor $d>1$ of $[G:S]$ with
$d\equiv1\pmod s$ satisfies $d\ge s+1$.

Condition~(i) is incompatible with Nakajima's bound unless
$\gamma(\cX)=0$ (Lemma~\ref{lem:positive}). For curves of $p$-rank
zero, the number of Sylow $p$-subgroups of $G$ is congruent to $1$
modulo $s$ (Lemma~\ref{lem:zero}); since it divides $[G:S]$,
condition~(ii) forces it to be $1$. Thus, under (i), $G$ fixes a
point of $\cX$ if and only if $S$ is normal in $G$, and (ii) is a
numerical condition guaranteeing normality. Condition~(i) cannot be
relaxed to $s\ge p(g-1)/(p-2)$: equality is attained by ordinary
curves of genus $p-1$ carrying a cyclic group of order $2p$ without
fixed points (Example~\ref{ex:sharp}). Condition~(ii) cannot be
omitted, as the Hermitian curve shows (Example~\ref{ex:hermitian}).
If one is willing to use \cite{GK10}, condition~(i) alone forces $G$
to fix a point unless $\cX$ is one of three classical curves
(Proposition~\ref{prop:GK}). Theorem~A replaces this identification of
the exceptions by a condition on $|G|$ that can be checked without
knowing $\cX$.

Our application concerns the generalized Suzuki curve. For an
integer $\ell\ge1$ put
\begin{equation}\label{eq:parameters}
 q_0=p^\ell,\qquad q=p^{2\ell-1},
\end{equation}
and let $\cX_\GS$ be the smooth projective model over $\F_q$ of
\begin{equation}\label{eq:GSmodel}
 y^q-y=x^{q_0}(x^q-x).
\end{equation}
Its genus is $q_0(q-1)/2$. The automorphisms in \eqref{eq:GSgroup}
below form a group $\Gamma(q,q_0)$ of order $q^2(q-1)$, and for
$\ell\ge2$ this is the full geometric automorphism group of
$\cX_\GS$ \cite[Theorem~1.3]{BC}; for $\ell=1$ the automorphism group
is larger (Remark~\ref{rem:l1}).

\begin{theoremB}
Let $\cX$ be a curve over $K$ of genus $q_0(q-1)/2$, with $q_0,q$ as
in \eqref{eq:parameters}, and let $G\le\Aut(\cX)$ have order
$q^2(q-1)$. Then $\gamma(\cX)=0$, the Sylow $p$-subgroup $S$ of $G$
is normal of order $q^2$, $G=S\rtimes C$ with $C$ cyclic of order
$q-1$, and $G$ fixes exactly one point $P$ of $\cX$. If
$K=\overline{\F}_q$, the curve $\cX$ is defined over $\F_q$, and $G$
is stable under the $q$-Frobenius, for instance if
$G=\Aut_{\overline{\F}_q}(\cX)$ or if every element of $G$ is defined
over $\F_q$, then $P\in\cX(\F_q)$.
\end{theoremB}

Using the fixed point provided by Theorem~B and the ramification
filtration at that point, we recover the curve itself.

\begin{theoremC}
Let $\ell\ge2$, let $q_0,q$ be as in \eqref{eq:parameters}, and let
$\cX$ be a curve over $K$ of genus $q_0(q-1)/2$. Then $\Aut(\cX)$
contains a subgroup of order $q^2(q-1)$ if and only if $\cX$ is
isomorphic to $\cX_\GS\otimes K$. In this case
$\Aut(\cX)\cong\Gamma(q,q_0)$.
\end{theoremC}

\begin{corollaryD}
Let $\ell\ge2$, let $q_0,q$ be as in \eqref{eq:parameters}, and let
$\cX/\F_q$ be a smooth projective geometrically irreducible curve of
genus $q_0(q-1)/2$. Suppose that $\Aut_{\overline{\F}_q}(\cX)$
contains a subgroup of order $q^2(q-1)$ all of whose elements are
defined over $\F_q$. Then $\cX$ is $\F_q$-isomorphic to $\cX_\GS$.
Moreover, $\cX$ has $p$-rank zero, and $\Aut_{\overline{\F}_q}(\cX)$
fixes exactly one point of $\cX$; this point lies in $\cX(\F_q)$.
\end{corollaryD}

In \cite[Theorem~1.1]{Taf24}, the curve $\cX_\GS$ was characterized
among curves over $\F_q$ of genus $q_0(q-1)/2$ by the order of its
automorphism group, under the additional hypotheses that all
automorphisms are defined over $\F_q$ and that the automorphism group
fixes an $\F_q$-rational point. Corollary~D shows that the
fixed-point hypothesis is superfluous, and that it suffices to have
a subgroup of the given order consisting of $\F_q$-automorphisms. The
proof of Theorem~C given in Section~\ref{sec:char} is self-contained
and does not rely on \cite{Taf24}. The hypothesis that the
automorphisms are defined over $\F_q$ cannot be dropped: nontrivial
twists of $\cX_\GS$ satisfy all the other hypotheses of Corollary~D
and the conclusions of Theorem~B, but are not $\F_q$-isomorphic to
$\cX_\GS$ (Proposition~\ref{prop:twist}).

\section{Preliminaries}\label{sec:prelim}

Throughout, $K$ is an algebraically closed field of odd characteristic
$p$, and $\cX/K$ is a smooth projective curve of genus $g\ge2$.
We write $g(\mathcal Y)$ and $\gamma(\mathcal Y)$ for the genus and the $p$-rank
of a curve $\mathcal Y$, respectively. Automorphisms act on the function
field $K(\cX)$. For $\beta\in K$ and $r\ge0$, the notation
$\beta^{p^{-r}}$ denotes the unique $p^r$-th root of $\beta$ in $K$.

\subsection{Quotients, stabilizers, and $p$-rank}

For a finite group $A\le\Aut(\cX)$ we write $\cX/A$ for the quotient
curve and $A_P$ for the stabilizer of $P\in\cX$. The stabilizer $A_P$
has a unique Sylow $p$-subgroup, namely its first ramification group
$A_P^{(1)}$, and $A_P/A_P^{(1)}$ is cyclic of order prime to $p$
\cite[Chapter~IV, \S2]{Serre}. If $S$ is a $p$-group of order $s$
acting on $\cX$, with short orbits of lengths $\ell_1,\ldots,\ell_b$,
the Deuring--Shafarevich formula reads
\begin{equation}\label{eq:DS}
 \gamma(\cX)-1=s\bigl(\gamma(\cX/S)-1\bigr)+\sum_{i=1}^b(s-\ell_i);
\end{equation}
see \cite[Chapter~11]{HKT}. Each $\ell_i$ is a proper divisor of $s$,
so $\ell_i\le s/p$.

\subsection{Ramification groups and the different}

Let $A\le\Aut(\cX)$ be finite, let $P\in\cX$, and put $G=A_P$.
For a local parameter $\pi$ at $P$, write
$G_0=G\ge G_1\ge G_2\ge\cdots$ for the lower ramification groups,
where $G_i$ consists of the automorphisms $h\in G$ for which
$v_P(h(\pi)-\pi)\ge i+1$ (for $i\ge0$).
In particular, $G_1$ is the unique Sylow $p$-subgroup of $G$.
There are injective homomorphisms
\[
 \theta_0\colon G_0/G_1\longrightarrow K^\times,\qquad
 \theta_i\colon G_i/G_{i+1}\longrightarrow(K,+)\quad(i\ge1),
\]
where $\theta_0(h)$ is the residue at $P$ of $h(\pi)/\pi$.
For $i\ge1$, the map $\theta_i$ is induced by the coefficient of
$\pi^{i+1}$ in $h(\pi)-\pi$ (with the usual convention for the action
on local parameters). These maps satisfy
\begin{equation}\label{eq:conj}
 \theta_i(chc^{-1})=\theta_0(c)^i\,\theta_i(h)
 \qquad(c\in G,\ h\in G_i,\ i\ge1).
\end{equation}
In particular, $G_i/
G_{i+1}$ is an elementary abelian $p$-group
for $i\ge1$; see \cite[Chapter~IV, \S2]{Serre}.

For a subgroup $U\le G$, its lower ramification groups at $P$ are
$U\cap G_i$. Hilbert's different formula gives
\[
 d_P(U)=\sum_{i\ge0}(|U\cap G_i|-1),
\]
where $d_P(U)$ is the different exponent of the cover
$\cX\to\cX/U$ at $P$; see \cite[Chapter~IV, \S1]{Serre}.
If this cover is unramified outside $P$, the Riemann--Hurwitz formula
therefore becomes
\begin{equation}\label{eq:RH}
 2g-2=|U|\bigl(2g(\cX/U)-2\bigr)+d_P(U).
\end{equation}

\section{Fixed points and \texorpdfstring{$p$}{p}-rank}\label{sec:lemmas}

The following lemma is well known (it is used, for instance, in
\cite{GKzero}); we include the short proof for completeness.

\begin{lemma}\label{lem:zero}
Assume $\gamma(\cX)=0$.
\begin{enumerate}[label=\textup{(\alph*)},leftmargin=*]
\item Every nontrivial $p$-subgroup of $\Aut(\cX)$ fixes exactly one
point of $\cX$ and acts freely on its complement.
\item Let $G\le\Aut(\cX)$, $S\in\Syl_p(G)$ and $s=|S|>1$, and let $n$
be the number of Sylow $p$-subgroups of $G$. Then $n\equiv1\pmod s$.
If $S$ is normal in $G$, then the fixed point of $S$ is the unique
fixed point of $G$, and $G=S\rtimes C$ with $C$ cyclic of order
$[G:S]$. If $S$ is not normal in $G$, then $G$ has no fixed point.
\end{enumerate}
\end{lemma}

\begin{proof}
(a) Let $\sigma$ have order $p$, with $k$ fixed points, and put
$\gamma'=\gamma(\cX/\langle\sigma\rangle)$. By \eqref{eq:DS},
$-1=p(\gamma'-1)+k(p-1)$. If $\gamma'\ge1$, the right side is
nonnegative; hence $\gamma'=0$ and $k=1$. Let $T\ne1$ be a $p$-group
and choose $z\in Z(T)$ of order $p$, with fixed point $Q$. Since $T$
commutes with $z$, it permutes the fixed points of $z$, so $T$ fixes
$Q$. A nonidentity $t\in T$ has a power of order $p$; that power fixes
only $Q$, hence $t$ fixes no other point.

(b) For $T\in\Syl_p(G)$ let $P_T$ be its fixed point. If
$P_T=P_{T'}=P$, then $T$ and $T'$ are $p$-subgroups of $G_P$ of the
maximal order $s$, so both equal the unique Sylow $p$-subgroup of
$G_P$. Thus $T\mapsto P_T$ is injective, and its image $\Omega$ is a
$G$-orbit of length $n$, because $gP_T=P_{gTg^{-1}}$. By (a), $S$
fixes $P_S$ and acts freely on $\Omega\setminus\{P_S\}$, so
$s\mid n-1$.

If $S$ is normal, then $G$ fixes $P=P_S$, and any further fixed point
of $G$ would be a second fixed point of $S$. Moreover $G=G_P$, so $S$
is the first ramification group of $G_P$ and $G/S$ is cyclic of order
$m=[G:S]$, prime to $p$. Choose $g\in G$ with $gS$ generating $G/S$.
Since $g^m\in S$, the order of $g$ is $mp^j$ for some $j\ge0$, so
$h=g^{p^j}$ has order $m$, and $hS$ generates $G/S$ because
$\gcd(p^j,m)=1$. Hence $C=\langle h\rangle$ is a complement to $S$.

If $S$ is not normal, then $n>1$. A fixed point of $G$ would be fixed
by every Sylow $p$-subgroup, hence would equal every $P_T$, which is
impossible.
\end{proof}

Part (a) of the next lemma is Nakajima's bound \cite{Nak}; curves
attaining it are studied in \cite{GK17}. We include the proof, since
the same computation gives part (b), which we also need.

\begin{lemma}\label{lem:positive}
Let $S$ be a $p$-group of order $s>1$ acting on $\cX$, and assume
$\gamma(\cX)\ge1$.
\begin{enumerate}[label=\textup{(\alph*)},leftmargin=*]
\item If $\gamma(\cX)\ge2$, then $s\le p(\gamma(\cX)-1)/(p-2)$.
\item If $\gamma(\cX)=1$, then $S$ acts freely and $s\le g-1$.
\end{enumerate}
In either case $s\le p(g-1)/(p-2)$.
\end{lemma}

\begin{proof}
Put $\bar\gamma=\gamma(\cX/S)$ and use \eqref{eq:DS} with
$\ell_i\le s/p$. If $\bar\gamma\ge2$, then $\gamma(\cX)-1\ge s$. If
$\bar\gamma=1$ and $b\ge1$, then $\gamma(\cX)-1\ge s(1-1/p)$. If
$\bar\gamma=0$ and $b=0$, formula \eqref{eq:DS} would give
$\gamma(\cX)=1-s<0$, a contradiction. If $\bar\gamma=0$ and $b=1$,
it gives $\gamma(\cX)=1-\ell_1\le0$, contrary to $\gamma(\cX)\ge1$.
Hence $b\ge2$ when $\bar\gamma=0$, and then
$\gamma(\cX)-1\ge s(1-2/p)$. In each of these three cases
$\gamma(\cX)-1>0$, so $\gamma(\cX)\ge2$, and the weakest of the three
bounds gives (a).

The only remaining possibility is $\bar\gamma=1$ and $b=0$. Then
$\gamma(\cX)=1$, the cover $\cX\to\cX/S$ is \'etale, and the Hurwitz
formula gives $g-1=s\bigl(g(\cX/S)-1\bigr)$. Since $g\ge2$, we get
$g(\cX/S)\ge2$ and $s\le g-1$. This proves (b), and shows that
$\gamma(\cX)=1$ occurs only in this case. The final assertion follows
from $\gamma(\cX)\le g$ and $g-1<p(g-1)/(p-2)$.
\end{proof}

\section{The numerical criterion}\label{sec:proofA}

\begin{proof}[Proof of Theorem A]
By Lemma~\ref{lem:positive}, condition~(i) excludes
$\gamma(\cX)\ge1$, so $\gamma(\cX)=0$. Let $n$ be the number of Sylow
$p$-subgroups of $G$. By Lemma~\ref{lem:zero}(b), $n\equiv1\pmod s$.
On the other hand, $n=[G:N_G(S)]$ divides $[G:S]$. Condition~(ii)
therefore gives $n=1$, that is, $S$ is normal in $G$, and
Lemma~\ref{lem:zero}(b) yields the decomposition $G=S\rtimes C$ and
shows that the fixed point $P$ of $S$ is the unique fixed point of
$G$.

Assume now that $K=\overline{\F}$, that $\cX$ is defined over $\F$,
and let $\varphi$ be the Frobenius map of $\cX(\overline{\F})$
relative to $\F$. By hypothesis $\varphi G\varphi^{-1}=G$. For
$\tau\in G$ we have
$\tau(\varphi(P))=\varphi\bigl((\varphi^{-1}\tau\varphi)(P)\bigr)
=\varphi(P)$, because $\varphi^{-1}\tau\varphi\in G$. Thus
$\varphi(P)$ is a fixed point of $G$, and uniqueness gives
$\varphi(P)=P$, that is, $P\in\cX(\F)$.
\end{proof}

\begin{remark}\label{rem:ii}
Condition~(ii) is strictly weaker than $|G|<s(s+1)$. For instance, if
$s=9$ and $[G:S]=16$, then $|G|=144>s(s+1)$, but no divisor $d>1$ of
$16$ is congruent to $1$ modulo $9$.
\end{remark}

\begin{example}[Condition (i) cannot be relaxed]\label{ex:sharp}
Let $\cX$ be the smooth projective model of
\[
 y^p-y=x+\frac1x .
\]
The right side has simple poles at $x=0$ and $x=\infty$ and no other
poles. Hence $\cX\to\PP^1_x$ is an Artin--Schreier cover of degree $p$,
totally ramified exactly over $0$ and $\infty$, and
\[
 g=\frac{p-1}{2}\bigl(-2+(1+1)+(1+1)\bigr)=p-1
\]
by \cite[Proposition~3.7.8]{Sti}. For $S=\langle\sigma\rangle$ with
$\sigma(y)=y+1$, formula \eqref{eq:DS} gives
$\gamma(\cX)-1=-p+2(p-1)$, so $\gamma(\cX)=p-1=g$. The involution
$\iota\colon(x,y)\mapsto(1/x,y)$ commutes with $\sigma$, and
$G=\langle\sigma,\iota\rangle$ is cyclic of order $2p$. Here $s=p$,
$[G:S]=2$ and
\[
 \frac{p}{p-2}(g-1)=p=s.
\]
Thus (ii) holds, since $2\not\equiv1\pmod p$, and (i) fails only by
equality, but $\cX$ is ordinary and $G$ has no fixed point: the two
fixed points of $S$, above $x=0$ and $x=\infty$, are interchanged by
$\iota$. Note that $S$ is normal in $G$; the failure is entirely due
to $\gamma(\cX)>0$. This curve attains Nakajima's bound; it falls
under case~(ii) of \cite[Theorem~1.1]{GK17}, and the family
$y^p-y=ax+1/x$ is discussed in \cite[Remark~2.2]{GK17}.
\end{example}

\begin{example}[Condition (ii) cannot be omitted]\label{ex:hermitian}
Let $\cX$ be the Hermitian curve $y^n+y=x^{n+1}$ with $n=p^h\ge3$,
of genus $n(n-1)/2$ and $p$-rank zero, and let $G=\Aut(\cX)$, which
is isomorphic to $\PGU(3,n)$; see \cite[Chapter~12]{HKT}. Then
$s=n^3$ satisfies~(i), since
$p(g-1)/(p-2)\le3(g-1)<3n^2/2<n^3$. However, $G$ acts transitively
on the $n^3+1$ points of $\cX(\F_{n^2})$ and has no fixed point;
accordingly (ii) fails: $d=n^3+1$ divides
$[G:S]=(n^3+1)(n^2-1)$ and $d\equiv1\pmod s$.
\end{example}

The next proposition shows what condition~(i) yields if the
classification of \cite{GK10} is used.

\begin{proposition}\label{prop:GK}
Let $G\le\Aut(\cX)$, let $S\in\Syl_p(G)$, and suppose that $s=|S|$
satisfies condition~\textup{(i)} of Theorem~A. If $G$ fixes no point
of $\cX$, then $\cX$ is birationally equivalent to one of the
following curves:
\begin{enumerate}[label=\textup{(\alph*)},leftmargin=*]
\item for $p\ge5$, the hyperelliptic curve $y^2=x^n-x$, where $n$ is
a power of $p$;
\item the Hermitian curve $y^n+y=x^{n+1}$, where $n\ge3$ is a power
of $p$;
\item for $p=3$, the Deligne--Lusztig curve of Ree type.
\end{enumerate}
\end{proposition}

\begin{proof}
We first show that $s>pg/(p-1)$. If $g\ge p-1$, this follows from~(i),
because $p(g-1)/(p-2)\ge pg/(p-1)$ is equivalent to $g\ge p-1$. If
$g\le p-2$, then $pg/(p-1)<p\le s$.

By Lemma~\ref{lem:positive}, condition~(i) gives $\gamma(\cX)=0$, so
$S$ fixes a point $P$ by Lemma~\ref{lem:zero}(a). Thus the stabilizer
of $P$ in $\Aut(\cX)$ contains a $p$-subgroup of order greater than
$pg/(p-1)$. If $G$ fixes no point of $\cX$, then neither does
$\Aut(\cX)$. By the main theorem of \cite{GK10}, $\cX$ is then
birationally equivalent to one of the curves (II)--(V) listed there.
Curve (IV) is the Deligne--Lusztig curve of Suzuki type, which occurs
only for $p=2$; the remaining ones are the curves in (a)--(c), except
that (a) is listed there for every odd $p$.

It remains to exclude $p=3$ in (a). In this case $g=(n-1)/2\ge2$, so
$n\ge9$. The hyperelliptic involution $\iota$ is central in
$\Aut(\cX)$, and since $p$ is odd, $S$ embeds into
$\Aut(\cX)/\langle\iota\rangle\le\PGL(2,K)$ and permutes the $n+1$
branch points of $\cX\to\cX/\langle\iota\rangle\cong\PP^1$. A
nontrivial $p$-subgroup of $\PGL(2,K)$ is conjugate to a group of
translations $x\mapsto x+b$, so it fixes exactly one point of $\PP^1$
and acts freely on the others. As $p\nmid n+1$, the fixed point of $S$
is a branch point, and $S$ acts freely on the remaining $n$ branch
points; hence $s\le n$. But
$p(g-1)/(p-2)=3(n-3)/2\ge n$ for $n\ge9$, contrary to~(i).
\end{proof}

\begin{remark}
Proposition~\ref{prop:GK} identifies the exceptions to Theorem~A
without using condition~(ii). Its proof, however, depends on
\cite{GK10}, which in turn uses the Kantor--O'Nan--Seitz
classification \cite{KOS}. Theorem~A is weaker in this respect, but
self-contained. Moreover, condition~(ii) can be checked from $s$ and
$|G|$ alone, without any knowledge of $\cX$; this is the situation in
Section~\ref{sec:suzuki}.
\end{remark}

\section{The generalized Suzuki curve and Theorem B}\label{sec:suzuki}

Keep the notation \eqref{eq:parameters} and \eqref{eq:GSmodel}.
The smooth projective model of \eqref{eq:GSmodel} has genus
$q_0(q-1)/2$ (see \cite{BC}). Its automorphism group contains the
transformations
\begin{equation}\label{eq:GSgroup}
 (x,y)\longmapsto
 \bigl(ax+b,\;ab^{q_0}x+a^{q_0+1}y+c\bigr),
 \qquad a\in\F_q^\times,\quad b,c\in\F_q.
\end{equation}
Indeed, if $(X,Y)$ denotes the transformed coordinates, then
$X^q-X=a(x^q-x)$ and
\[
Y^q-Y=(ab^{q_0}+a^{q_0+1}x^{q_0})(x^q-x)
=X^{q_0}(X^q-X).
\]
These transformations form a group $\Gamma(q,q_0)$ of order
$q^2(q-1)$, whose Sylow $p$-subgroup (given by $a=1$) is normal of
order $q^2$. For $\ell\ge2$ it is the full geometric automorphism
group of $\cX_\GS$ \cite[Theorem~1.3]{BC}. For Theorem~B we need only
its order.

\begin{remark}\label{rem:l1}
The restriction $\ell\ge2$ is needed for the last statement. For
$\ell=1$ we have $q_0=q=p$, and in terms of $y_1=y-x^2$ equation
\eqref{eq:GSmodel} becomes $y_1^p-y_1=x^2-x^{p+1}$. Let $c$ be a root
of the separable additive polynomial $E(T)=T^{p^2}-2T^p+T$, put
$\beta=c^p-2c$, so that $\beta^p=-c$, and let $e$ satisfy
$e^p-e=c^2-c^{p+1}$. A direct computation shows that
$(x,y_1)\mapsto(x+c,\;y_1+\beta x+e)$ is an automorphism of this
curve. These maps form a $p$-group of order $p^3>q^2$, so
$\Aut(\cX_\GS)$ is strictly larger than $\Gamma(p,p)$ when $\ell=1$.
Theorem~B still applies to $\Gamma(p,p)$ and to every subgroup of
$\Aut(\cX_\GS)$ of order $p^2(p-1)$.
\end{remark}

\begin{proof}[Proof of Theorem B]
Since $p\nmid q-1$, we have $s=q^2$ and $[G:S]=q-1<s+1$, so (ii)
holds. Let $g=q_0(q-1)/2$. For (i):
\begin{itemize}[leftmargin=*]
\item If $p\ge5$, then $q_0\le q$, so $g-1<q_0q/2\le q^2/2$ and
$\frac{p}{p-2}(g-1)<\frac{p}{2(p-2)}\,q^2<q^2$.
\item If $p=3$ and $\ell\ge2$, then $q_0\le q/3$, so
$3(g-1)<\frac32q_0(q-1)\le\frac12q(q-1)<q^2$.
\item If $(p,\ell)=(3,1)$, then $q_0=q=3$, $g=3$ and
$3(g-1)=6<9=q^2$.
\end{itemize}
Theorem~A therefore applies, and gives $G=S\rtimes C$ with $C$ cyclic
of order $q-1$. If $\cX$ is defined over $\F_q$, then
$\Aut_{\overline{\F}_q}(\cX)$ is stable under the Frobenius; if every
element of $G$ is defined over $\F_q$, each of them commutes with the
Frobenius. In both cases $P\in\cX(\F_q)$ by Theorem~A.
\end{proof}

\begin{remark}\label{rem:GK19}
Theorem~B is not a consequence of \cite{GK19}. Since $q_0^2=pq$, we
have $g^2=pq(q-1)^2/4$, hence
\[
 |G|=q^2(q-1)=\frac{4q}{p(q-1)}\,g^2\le 2g^2,
\]
far below the threshold $252g^2$ of \cite{GK19}; moreover, no
solvability of $G$ is assumed.
\end{remark}

For $\cX=\cX_\GS$ and $\ell\ge2$, Theorem~B recovers the facts that
$\cX_\GS$ has $p$-rank zero and that its full automorphism group fixes
exactly one point, the point at infinity. The content of Theorem~B is
that these conclusions hold for every curve of the given genus and
every automorphism group of the given order.

\section{Proof of Theorem C}\label{sec:char}

The ``if'' part is clear, and the last assertion then follows from
\cite[Theorem~1.3]{BC}. For the converse, fix $\ell\ge2$, write
$n=2\ell-1$, so that $q=p^n$, and let $\cX$ be a curve over $K$ of
genus $g=q_0(q-1)/2$ with a subgroup $G\le\Aut(\cX)$ of order
$q^2(q-1)$. By Theorem~B, $\gamma(\cX)=0$, $G=S\rtimes C$ with
$|S|=q^2$ and $C$ cyclic of order $q-1$, and $G$ fixes exactly one
point $P$. By Lemma~\ref{lem:zero}(a), for every subgroup $1\ne U\le S$
the cover $\cX\to\cX/U$ is totally ramified at $P$ and unramified
elsewhere. We now apply the ramification notation from
Section~\ref{sec:prelim} to $G=G_P$. Write
$G_0=G\ge G_1=S\ge G_2\ge\cdots$ for the lower ramification
groups at $P$ and use $\theta_i$ for the associated homomorphisms.
Since $|G_0/G_1|=q-1$, the image of $\theta_0$ is $\F_q^\times$;
fix a generator $\zeta$ of this group. For each $1\ne U\le S$,
the cover $\cX\to\cX/U$ is ramified only at $P$, so
\eqref{eq:RH} applies.

\begin{lemma}\label{lem:filtration}
\begin{enumerate}[label=\textup{(\alph*)},leftmargin=*]
\item $g(\cX/S)=0$ and $\sum_{i\ge2}(|G_i|-1)=q_0(q-1)$.
\item $|G_2|=q$, $G_{q_0+1}=G_2$ and $G_{q_0+2}=1$.
\item $T:=G_2$ is an elementary abelian normal subgroup of $G$
contained in the center of $S$, the curve $\cX/T$ has genus $0$, and
$\theta_{q_0+1}(T)=\F_qw_0$ for some $w_0\in K^\times$.
\end{enumerate}
\end{lemma}

\begin{proof}
(a) The terms $i=0,1$ of $d_P(S)$ contribute $2(q^2-1)$. If
$g(\cX/S)\ge1$, then \eqref{eq:RH} gives $g-1\ge q^2-1$, whereas
$g<q^2$. Hence $g(\cX/S)=0$, so $d_P(S)=2g-2+2q^2$, and subtracting
the terms $i=0,1$ gives the second assertion.

(b) If $G_2=G_1$, the sum in (a) is at least $q^2-1>q_0(q-1)$.
Hence $V_1=\theta_1(G_1/G_2)$ is a nonzero subgroup of $K$, and by
\eqref{eq:conj} it is closed under multiplication by $\zeta$, so it is
an $\F_q$-subspace of $K$. Thus $|G_1/G_2|\in\{q,q^2\}$, and
$|G_1/G_2|=q^2$ would give $G_2=1$, contradicting (a). So $|G_2|=q$.

Let $m\ge2$ be maximal with $G_m=G_2$. By (a),
$(m-1)(q-1)\le q_0(q-1)$, so $m\le q_0+1$. By \eqref{eq:conj},
$V_m=\theta_m(G_m/G_{m+1})$ is a nonzero vector space over the field
$F'=\F_p(\zeta^m)\subseteq\F_q$. Suppose that $F'\ne\F_q$, so that
$F'=\F_{p^e}$ for a proper divisor $e$ of $n$. Then
$\zeta^{m(p^e-1)}=1$, so $q-1$ divides $m(p^e-1)$, that is,
$N=(q-1)/(p^e-1)$ divides $m$. Since $n$ is odd, $n/e\ge3$ and
$N=\sum_{j=0}^{n/e-1}p^{je}\ge1+p^e+p^{n-e}$; moreover
$e\le n/3\le\ell-1$, so $N\ge1+p+p^\ell>q_0+1\ge m$, contradicting
$N\mid m$. Hence $F'=\F_q$, so $|V_m|\ge q=|G_m|$, which forces
$G_{m+1}=1$. Now (a) gives $(m-1)(q-1)=q_0(q-1)$, that is,
$m=q_0+1$.

(c) By (b), $\theta_{q_0+1}$ is injective on $T=G_{q_0+1}$, so $T$ is
elementary abelian, and its image is an $\F_q$-subspace of $K$ of
order $q$, that is, $\F_qw_0$ for some $w_0\ne0$. Ramification groups
are normal in $G_0=G$. For $s\in S=G_1$ we have $\theta_0(s)=1$, so
\eqref{eq:conj} gives $\theta_{q_0+1}(sts^{-1})=\theta_{q_0+1}(t)$,
hence $sts^{-1}=t$ for $t\in T$. Finally, $T\cap G_i=T$ for
$0\le i\le q_0+1$ and $T\cap G_i=1$ otherwise, so
$d_P(T)=(q_0+2)(q-1)$, and \eqref{eq:RH} gives
$q_0(q-1)-2=2q\,g(\cX/T)-2q+(q_0+2)(q-1)$, that is, $g(\cX/T)=0$.
\end{proof}

Let $\lambda\colon T\to\F_q$, $\lambda(t)=w_0^{-1}\theta_{q_0+1}(t)$.
This is a group isomorphism, and by \eqref{eq:conj}
\begin{equation}\label{eq:lambda}
 \lambda(ctc^{-1})=\theta_0(c)^{q_0+1}\lambda(t)\qquad
 (c\in G,\ t\in T).
\end{equation}

\begin{lemma}\label{lem:x}
There is $x\in K(\cX)$ with $K(\cX)^T=K(x)$, whose only pole is $P$,
of order $q$, such that $c(x)=\theta_0(c)^{-1}x$ for $c\in C$ and
$\{s(x)-x:s\in S\}=\F_q$.
\end{lemma}

\begin{proof}
By Lemma~\ref{lem:filtration}(c), $\cX/T\cong\PP^1$, so
$K(\cX)^T=K(x_0)$ with $x_0$ having a single simple pole at the image
of $P$; as $P$ is totally ramified over it, $x_0$ has a pole of order
$q$ at $P$ and no other pole. The group $G/T$ acts faithfully on
$K(x_0)$ and fixes this pole, so each of its elements maps $x_0$ to
$ax_0+b$ with $a\in K^\times$, $b\in K$. Elements of $S/T$ have
$p$-power order, so for them $a=1$; as $|S/T|=q$, the translation
parts form an additive subgroup $W\subset K$ of order $q$. An element
of $C$ with $a=1$ would be a translation of order prime to $p$, hence
trivial; so a generator of $C$ has $a\ne1$ and fixes a point of the
affine line, and after translating $x_0$ we get $c(x_0)=a_cx_0$ for
all $c\in C$. Writing $x_0=u\pi^{-q}$ with $u$ a unit at $P$, we have
$c(x_0)/x_0=(c(u)/u)\,(c(\pi)/\pi)^{-q}$; since $c$ fixes $P$, the
first factor has value $1$ at $P$, so
$a_c=\theta_0(c)^{-q}=\theta_0(c)^{-1}$, because
$\theta_0(c)\in\F_q$. In particular $\{a_c:c\in C\}=\F_q^\times$. For
$s\in S$ with $s(x_0)=x_0+b$ one computes
$csc^{-1}(x_0)=x_0+a_c^{-1}b$, so $W$ is an $\F_q$-subspace of $K$
of order $q$, that is, $W=\F_qb_1$ with $b_1\ne0$. Then $x=x_0/b_1$
has the required properties.
\end{proof}

\begin{lemma}\label{lem:y}
There is $y\in K(\cX)$, regular outside $P$, with $K(\cX)=K(x,y)$,
such that
\begin{enumerate}[label=\textup{(\alph*)},leftmargin=*]
\item $t(y)=y+\lambda(t)$ for $t\in T$, and
$c(y)=\theta_0(c)^{-(q_0+1)}y$ for $c\in C$;
\item $y^q-y=f(x)$ with $0\ne f\in K[x]$, where
$f(\alpha x)=\alpha^{q_0+1}f(x)$ for all $\alpha\in\F_q^\times$, and
no exponent of a monomial of $f$ is divisible by $q$.
\end{enumerate}
\end{lemma}

\begin{proof}
Let $R$ be the ring of functions on $\cX$ regular outside $P$. Then
$R\cap K(x)=K[x]$. The ring $R$ is a finitely generated torsion-free,
hence free, $K[x]$-module. Since $\cX\to\cX/T$ is unramified above
the affine line, for each $a\in K$ the fibre over $x=a$ consists of
$q$ points $Q_1,\dots,Q_q$, and evaluation gives
$R/(x-a)R\cong K^q$. The $K[x]$-linear map
$\tr_T\colon R\to K[x]$, $r\mapsto\sum_{t\in T}t(r)$, satisfies
$\tr_T(r)(a)=\sum_jr(Q_j)$, so it is surjective modulo every $x-a$. 
 Since $\operatorname{tr}_T(R)$ is an ideal of $K[x]$, we may write $\operatorname{tr}_T(R)=(h(x))$. The surjectivity modulo every $x-a$ implies that $h(a)\neq 0$ for all $a\in K$. As $K$ is algebraically closed, $h(x)$ is a nonzero constant. Hence
$$
\operatorname{tr}_T(R)=K[x].
$$

Choose $u\in R$ with $\tr_T(u)=1$ and
put $y_0=-\sum_{s\in T}\lambda(s)\,s(u)\in R$. For $t\in T$,
\[
 t(y_0)=-\sum_{s\in T}\lambda(t^{-1}s)\,s(u)
 =y_0+\lambda(t)\tr_T(u)=y_0+\lambda(t).
\]
For $c\in C$ and $t\in T$, \eqref{eq:lambda} gives
$t(c(y_0))=c\bigl((c^{-1}tc)(y_0)\bigr)
=c(y_0)+\theta_0(c)^{-(q_0+1)}\lambda(t)$. Hence
\[
 y=\frac1{q-1}\sum_{c\in C}\theta_0(c)^{q_0+1}\,c(y_0)\in R
\]
satisfies $t(y)=y+\lambda(t)$, and a reindexing of the sum gives
$c(y)=\theta_0(c)^{-(q_0+1)}y$. This proves (a).

Since $\prod_{\beta\in\F_q}(Y+\beta)=Y^q-Y$, the element
$f=y^q-y=\prod_{t\in T}t(y)$ is $T$-invariant and lies in
$R\cap K(x)=K[x]$. As $y$ has $q$ distinct conjugates under $T$, we
have $K(\cX)=K(x,y)$; and $f\ne0$, since otherwise $y\in\F_q$. For
$c\in C$ put $\alpha=\theta_0(c)^{-1}\in\F_q^\times$. Then
$c(f)=f(\alpha x)$ by Lemma~\ref{lem:x}, and
$c(f)=c(y)^q-c(y)=\alpha^{q_0+1}f$; as $c$ ranges over $C$, $\alpha$
ranges over $\F_q^\times$. Consequently every exponent $j$ of $f$
satisfies $j\equiv q_0+1\pmod{q-1}$.

Finally, let $j$ be the largest exponent of $f$ divisible by $q$, with
coefficient $a_j$, and let $\mu\in K$ satisfy $\mu^q=-a_j$. Since
$j/q\equiv j\equiv q_0+1\pmod{q-1}$, replacing $y$ by
$y+\mu x^{j/q}$ preserves (a), and it replaces $f$ by
$f-a_jx^j-\mu x^{j/q}$. This removes the exponent $j$ and only adds
the smaller exponent $j/q$. Repeating this step, we obtain (b).
\end{proof}

\begin{lemma}\label{lem:shape}
$f=Ax^{q+q_0}+Bx^{q_0+1}$ for some $A,B\in K$.
\end{lemma}

\begin{proof}
Write $f=\sum_ja_jx^j$. For $\mu\in\F_q^\times$ put
$w_\mu=\sum_{r=0}^{n-1}(\mu y)^{p^r}$. Then
$w_\mu^p-w_\mu=\mu^qy^q-\mu y=\mu f$, and
$t(w_\mu)=w_\mu+\Tr_{\F_q/\F_p}(\mu\lambda(t))$ for $t\in T$. Let
$T_\mu$ be the kernel of $t\mapsto\Tr_{\F_q/\F_p}(\mu\lambda(t))$, a
subgroup of index $p$, and let $E_\mu=K(\cX)^{T_\mu}=K(x,w_\mu)$, of
degree $p$ over $K(x)$. Let $P_\mu$ be the place of $E_\mu$ below $P$,
and let $\delta_\mu$ be the different exponent of $E_\mu/K(x)$ at
$P_\mu$. As in Lemma~\ref{lem:filtration}(c), the different exponents
of $K(\cX)/K(x)$ and of $K(\cX)/E_\mu$ at $P$ are $(q_0+2)(q-1)$ and
$(q_0+2)(q/p-1)$, and the ramification index of $P$ over $P_\mu$ is
$q/p$. By transitivity of the different \cite[Section~3.4]{Sti},
\[
 (q_0+2)(q-1)=(q_0+2)(q/p-1)+(q/p)\,\delta_\mu,
 \qquad\text{so}\qquad \delta_\mu=(q_0+2)(p-1).
\]

Now let $\wp(h)=h^p-h$. Every exponent of $f$ can be written as
$j=ip^r$ with $p\nmid i$ and $0\le r\le n-1$, since $q\nmid j$. As
$\beta x^{ip^r}-\beta^{p^{-r}}x^i\in\wp(K[x])$, we have
$\mu f\equiv g_\mu\pmod{\wp(K[x])}$, where
\[
 g_\mu=\sum_{p\nmid i}c_i(\mu)\,x^i,\qquad
 c_i(\mu)=\sum_{r=0}^{n-1}\mu^{p^{-r}}a_{ip^r}^{p^{-r}} .
\]
Thus $E_\mu=K(x,w')$ with $w'^p-w'=g_\mu$. The place $P_\mu$ is
ramified over $K(x)$, so $g_\mu$ is not constant; hence
$k_\mu=\deg g_\mu$ is prime to $p$, and \cite[Proposition~3.7.8]{Sti}
gives $\delta_\mu=(p-1)(k_\mu+1)$. Therefore $k_\mu=q_0+1$ for every
$\mu\in\F_q^\times$, so $c_i(\mu)=0$ for all $\mu\in\F_q$ whenever
$i>q_0+1$ and $p\nmid i$. The maps $\mu\mapsto\mu^{p^{-r}}$,
$0\le r\le n-1$, are distinct field embeddings $\F_q\to K$, hence
linearly independent over $K$ by Dedekind's lemma. It follows that
$a_{ip^r}=0$ for all such $i$ and all $r$.

So every exponent of $f$ has the form $j=ip^r$ with $p\nmid i$,
$i\le q_0+1$ and $0\le r\le n-1$, and $j\equiv q_0+1\pmod{q-1}$ by
Lemma~\ref{lem:y}. Multiplying by $p^{n-r}$ and using $p^n\equiv1$, we
get $i\equiv p^{\bar\alpha}+p^{\bar\beta}\pmod{q-1}$, where
$\bar\alpha,\bar\beta\in\{0,\dots,n-1\}$ are the residues of $\ell-r$
and $-r$ modulo $n$. Both sides lie in $[1,q-2]$, because
$q_0+1<q-1$ and $2p^{n-1}<q-1$; hence $i=p^{\bar\alpha}+p^{\bar\beta}$.
As $p\nmid i$ and $\ell\not\equiv0\pmod n$, exactly one of
$\bar\alpha,\bar\beta$ vanishes. If $\bar\beta=0$, then $r=0$ and
$i=q_0+1$, so $j=q_0+1$. If $\bar\alpha=0$, then $r=\ell$,
$\bar\beta=\ell-1$ and $i=1+p^{\ell-1}$, so $j=p^\ell+p^{2\ell-1}=q+q_0$.
\end{proof}

\begin{proof}[Proof of Theorem C]
By Lemma~\ref{lem:shape}, $f=Ax^{q+q_0}+Bx^{q_0+1}$. Choose $s\in S$
with $s(x)=x+b$ and $b\in\F_q^\times$ (Lemma~\ref{lem:x}). By
Lemma~\ref{lem:filtration}(c), $s$ commutes with every $t\in T$, so
$t(s(y))=s(y)+\lambda(t)$; hence $h=s(y)-y$ is $T$-invariant and
regular outside $P$, that is, $h\in K[x]$. Applying $s$ to
$y^q-y=f(x)$ gives $h^q-h=f(x+b)-f(x)$. Since $b^q=b$, we have
$(x+b)^{q+q_0}=(x^q+b)(x^{q_0}+b^{q_0})$, and therefore
\[
 f(x+b)-f(x)=Ab^{q_0}x^q+(A+B)\,b\,x^{q_0}+Bb^{q_0}x+(A+B)\,b^{q_0+1}.
\]
Let $\nu\in K$ with $\nu^q=Ab^{q_0}$. Then
\[
 (h-\nu x)^q-(h-\nu x)=(A+B)\,b\,x^{q_0}+(Bb^{q_0}+\nu)\,x+(A+B)\,b^{q_0+1}.
\]
The right side has degree less than $q$, so $h-\nu x$ is constant and
the right side is constant. Hence $A=-B$ and $\nu=-Bb^{q_0}$. Raising
the latter to the $q$-th power and using $b^{q_0q}=b^{q_0}$ gives
$-Bb^{q_0}=Ab^{q_0}=\nu^q=-B^qb^{q_0}$, so $B\in\F_q$; and $B\ne0$
because $f\ne0$. Put $y'=-B^{-1}y$. Since $-B^{-1}\in\F_q$,
\[
 y'^q-y'=-B^{-1}f(x)=x^{q+q_0}-x^{q_0+1}=x^{q_0}(x^q-x).
\]
As $K(\cX)=K(x,y')$ has degree $q$ over $K(x)$, this is an
isomorphism of $K(\cX)$ with the function field of
$\cX_\GS\otimes K$, and so $\cX\cong\cX_\GS\otimes K$.
\end{proof}

\section{Curves over \texorpdfstring{$\F_q$}{Fq}}\label{sec:Fq}

\begin{lemma}\label{lem:center}
For $\ell\ge2$, the group $\Gamma=\Gamma(q,q_0)$ has trivial center.
\end{lemma}

\begin{proof}
Write $g_{a,b,c}$ for the map \eqref{eq:GSgroup}, and let
$z=g_{a,b,c}$ be central in $\Gamma$. Comparing $z\circ g_{1,\beta,0}$
with $g_{1,\beta,0}\circ z$ in the first coordinate gives $a\beta=\beta$
for all $\beta\in\F_q$, so $a=1$; in the second coordinate it then gives
$b^{q_0}\beta=\beta^{q_0}b$ for all $\beta\in\F_q$. If $b\ne0$, then
$r^{q_0}=r$ for all $r\in\F_q$, that is, $\F_q\subseteq\F_{p^\ell}$,
which is impossible because $2\ell-1>\ell$. Hence $b=0$. Comparing
$z=g_{1,0,c}$ with $g_{\alpha,0,0}$ gives $\alpha^{q_0+1}c=c$ for all
$\alpha\in\F_q^\times$. Since $0<q_0+1<q-1$, a generator $\alpha$ of
$\F_q^\times$ satisfies $\alpha^{q_0+1}\ne1$, so $c=0$ and $z=1$.
\end{proof}

Let $\sigma$ be the $q$-Frobenius in
$\operatorname{Gal}(\overline{\F}_q/\F_q)$. For a morphism $\eta$
between curves defined over $\F_q$ we write $\eta^\sigma$ for its
conjugate; $\eta$ is defined over $\F_q$ if and only if
$\eta^\sigma=\eta$.

\begin{proof}[Proof of Corollary D]
By Theorem~C there is an $\overline{\F}_q$-isomorphism
$\eta\colon\cX\to\cX_\GS$, and $\Aut_{\overline{\F}_q}(\cX)$ has
order $q^2(q-1)$. Hence the given subgroup is the whole automorphism
group, and every automorphism of $\cX$ is defined over $\F_q$. Every
element of $\Gamma=\Aut_{\overline{\F}_q}(\cX_\GS)$ is defined over
$\F_q$ \cite[Theorem~1.3]{BC}. Put $\tau=\eta^\sigma\eta^{-1}\in\Gamma$.
For $\alpha\in\Aut_{\overline{\F}_q}(\cX)$ let
$\alpha'=\eta\alpha\eta^{-1}\in\Gamma$. Since $(\alpha')^\sigma=\alpha'$,
\[
 \alpha^\sigma=(\eta^\sigma)^{-1}\alpha'\eta^\sigma
 =\eta^{-1}\,\tau^{-1}\alpha'\tau\,\eta .
\]
As $\alpha^\sigma=\alpha$ for every $\alpha$, the element $\tau$
commutes with every element of $\Gamma$, so $\tau=1$ by
Lemma~\ref{lem:center}. Thus $\eta^\sigma=\eta$, and $\eta$ is defined
over $\F_q$. The remaining assertions follow from Theorem~B.
\end{proof}

Frobenius-stability of $G$ suffices for the rationality of its unique
fixed point in Theorem~B, and it is weaker than requiring every
element of $G$ to be defined over $\F_q$. The next proposition shows
that the latter assumption cannot be removed from Corollary~D.

\begin{proposition}\label{prop:twist}
Let $\ell\ge2$ and $\Gamma=\Gamma(q,q_0)$. For every
$\tau\in\Gamma\setminus\{1\}$ there is a smooth projective
geometrically irreducible curve $\cX_\tau$ over $\F_q$ with the
following properties:
\begin{enumerate}[label=\textup{(\alph*)},leftmargin=*]
\item $\cX_\tau$ is $\overline{\F}_q$-isomorphic to $\cX_\GS$; in
particular it has genus $q_0(q-1)/2$ and
$\Aut_{\overline{\F}_q}(\cX_\tau)\cong\Gamma$;
\item the automorphisms of $\cX_\tau$ defined over $\F_q$ form a
subgroup isomorphic to the centralizer $C_\Gamma(\tau)$, which is a
proper subgroup of $\Gamma$;
\item $\cX_\tau$ is not $\F_q$-isomorphic to $\cX_\GS$.
\end{enumerate}
Since $\cX_\tau$ is defined over $\F_q$, its automorphism group is
Frobenius-stable, so Theorem~B applies to it. Thus, in Corollary~D,
the assumption that the automorphisms are defined over $\F_q$ cannot
be dropped.
\end{proposition}

\begin{proof}
Every element of $\Gamma$ is defined over $\F_q$, and
$\Gamma=\Aut_{\overline{\F}_q}(\cX_\GS)$ by \cite[Theorem~1.3]{BC},
so the Galois group acts trivially on $\Aut_{\overline{\F}_q}(\cX_\GS)$.
Hence the continuous homomorphism
$\operatorname{Gal}(\overline{\F}_q/\F_q)\to\Gamma$ sending $\sigma$ to
$\tau$ is a $1$-cocycle, and by \cite[Theorem~X.2.2]{Sil} there are a
curve $\cX_\tau$ over $\F_q$ and an $\overline{\F}_q$-isomorphism
$\eta\colon\cX_\tau\to\cX_\GS$ with $\eta^\sigma=\tau\circ\eta$. This
gives (a). For $\alpha\in\Aut_{\overline{\F}_q}(\cX_\tau)$ put
$\alpha'=\eta\alpha\eta^{-1}\in\Gamma$. As in the proof of
Corollary~D, $\alpha^\sigma=\eta^{-1}\tau^{-1}\alpha'\tau\eta$, so
$\alpha$ is defined over $\F_q$ if and only if
$\alpha'\in C_\Gamma(\tau)$. By Lemma~\ref{lem:center} and $\tau\ne1$,
this centralizer is proper, which proves (b). If $\cX_\tau$ were
$\F_q$-isomorphic to $\cX_\GS$, all its automorphisms would be defined
over $\F_q$, contrary to (b); this proves (c).
\end{proof}

\subsection*{Acknowledgements}
The first author was supported by CNPq grant
no.~302774/2025-4 and FAPESP grant no.~2024/00923-6.
The second author was supported by FAPESP grant
no.~2026/09571-0.

\textbf{Competing interests.}
The authors have no competing interests to declare.

\textbf{Data availability.} No data were generated or analysed.

\end{document}